\documentclass[12pt, twoside]{article}
\usepackage{graphicx}
\usepackage{hyperref}

\usepackage{graphicx}
\usepackage{hyperref}

\usepackage{graphicx}
\usepackage{hyperref}

\newtheorem{Theorem}{Theorem}

\newtheorem{Conjecture}{Conjecture}

\def\LAMBDA{\mbox{\rlap{$\raise3pt\hbox{--}$}{$\lambda$}}}  

\begin{document}

\begin{center}

{\Large\bf A Statistical Approach to Prime Gaps and Andrica's Conjecture} \\

\bigskip

Sameen Ahmed Khan\footnote{\normalsize
{\bf E-mail address:} \url{rohelakhan@yahoo.com} \\
{\bf URL:} \url{http://orcid.org/0000-0003-1264-2302}.  
} \\ 
Department of Mathematics and Sciences \\
College of Arts and Applied Sciences (CAAS) \\
Dhofar University \\
Post Box No. 2509, 
Postal Code: 211 \\
Salalah, Sultanate of Oman. \\

\end{center}

\medskip

\begin{abstract} 
We examine the prime gaps using a statistical approach.  
It is first shown that the Andrica's conjecture is true for half or more cases.  
Using the arguments of averages, it is further shown that Andrica's conjecture is true.  
We further obtain a precise bound for the Andrica's expression.  
\end{abstract}

\noindent
{\bf Keywords and phrases:} Prime Gap; Andrica's Conjecture. \\ 

\noindent
{\bf AMS Classification:} 11A41; 11N05. \\

\noindent
{\bf PACS:} 02.10.De 

\tableofcontents

\newpage

\setcounter{section}{0}

\section{Introduction} 
It is well known that there is no formula to find the $n$-th prime.  
So, the gaps between two successive primes are of a keen interest.  
We shall infer the average gap between the primes using statistical techniques.  
Here, we shall consider the Andrica's conjecture.  
Let the $n$-th prime be denoted by $p_n$.  
Andrica's conjecture states that 
\begin{Conjecture} 
{\bf (Andrica's Conjecture):} \label{theorem-andricas-conjecture}
Andrica's conjecture states that, for $p_n$ the $n$-th prime number, the inequality
$$
h_n \equiv \sqrt{p_{n + 1}} - \sqrt{p_n} < 1$$
holds for all $n$.  
The $\lim h_n = 0$. 
\end{Conjecture}
If $g_n = p_{n + 1} - p_n$ denotes the $n$-th prime gap, 
then Andrica's conjecture can also be rewritten as $g_n < 1 + 2 \sqrt{p_n}$. 
Numerical data confirms the conjecture for $n$ up to $1.3002 \times 10^{16}$ 
(see~\cite{Andrica-Wiki} for details).  
Dorin Andrica published this conjecture in 1986~\cite{Andrica-1986}. 
Andrica's conjecture implies other conjectures such as the Legendre's conjecture, which 
states that there is a prime number between $n^2$ and $(n+ 1)^2$ for every positive 
integer $n$.  It also has a bearing on certain other conjectures as well. 
In this article, we shall conclude that the average value of $\{ h_n \}$ 
denoted by $\bar{h}_n$ satisfies $\bar{h}_n \in (0\,, 1)$.  
From this, we further conclude that half or more of the $\{h_n \}$ are less than one. 
We work out an inequality for the average function, which states that the average value 
of $\{ h_n \}$ is a decreasing function.  
Using the arguments of averages, it is further shown that Andrica's conjecture is true.  
We further obtain a precise bound for the Andrica's expression.  
 
In the remainder of this section, we shall note a few basic results for the $k$-th prime 
and the approximations we may be using.  
In the next section, we shall present the statistical approach to the prime gaps.  
Section-3 has the statistical approach to the Andrica's conjecture.  
Section-4 has the proof of the Andrica's conjecture.  
Section-5 addresses the generalization of Andrica's conjecture. 
Section-5 has the concluding remarks.  
 
There are several bounds for the $k$-th prime such as $p_k > k \ln k$ due to 
Rosser~\cite{Rosser-1938}.  This was subsequently improved to~\cite{Dusart-MC}  
\begin{eqnarray}
k \left\{\ln k + \ln \ln k - 1 \right\}
< p_k < k \left\{\ln k + \ln \ln k \right\}\,, ~{\rm for}~ k \ge 6\,. 
\label{bounds-of-pk-big} 
\end{eqnarray} 
Better bounds are available in the works of Pierre Dusart~\cite{Dusart-RJ}  
\begin{eqnarray}
p_k 
& \le & 
k \left\{\ln k + \ln \ln k - 1 + \frac{\ln \ln k - 2}{\ln k}\right\} ~{\rm for}~ k \ge 688383 \nonumber \\
p_k 
& \ge & 
k \left\{\ln k + \ln \ln k - 1 + \frac{\ln \ln k - 2.1}{\ln k}\right\} ~{\rm for}~ k \ge 3\,.
\label{bounds-of-pk-huge} 
\end{eqnarray} 

In our approach, it is sufficient to use simpler expressions for the $p_k$.  
In some of the expressions, we shall use an ``oversimplification''.  
It is straightforward to see that 
\begin{eqnarray}
\left(\ln k + \ln \ln k \right) = \ln \left(k \ln k \right) < k\,.  
\label{upper-bound-of-right-bracket}
\end{eqnarray}
This simplification leads to 
\begin{eqnarray}
p_k < k^2\,.  
\label{upper-bound-of-pk-square}
\end{eqnarray}
The other simplification, we shall be using is 
\begin{eqnarray}
p_k \sim k \ln k\,.  
\label{pk-simple} 
\end{eqnarray}
%

 
\section{Statistical Approach to Prime Gaps} 
Let us consider the following sum 
\begin{eqnarray}
\sum_{k = 1}^{n} g_k 
& = & g_1 + g_2 + g_3 + \cdots + g_n \nonumber \\
& = & 
\left(p_2 - p_1 \right)
+ \left(p_3 - p_2 \right) 
+ \left(p_4 - p_3 \right)
+ \cdots \nonumber \\ 
& & \qquad \qquad \qquad 
+ \left(p_{n + 1} - p_n \right) \nonumber \\ 
& = & 
p_{n + 1} - p_1\,,  
\label{Sum-g}
\end{eqnarray}
where $p_1 = 2$.  We have used the telescoping property in the above summation.  
Next, we rewrite Eq.~(\ref{Sum-g}) as 
\begin{eqnarray}
\sum_{k = 1}^{n} g_k + 2 = p_{n + 1}\,.    
\label{Sum-g-t}
\end{eqnarray}
We divide the expression in Eq.~(\ref{Sum-g-t}) with $n$ and then rearrange to obtain 
the average value of the first $n$ prime gaps denoted by $\bar{g}_n$ 
\begin{eqnarray}
\bar{g}_n 
& \equiv & 
\frac{1}{n} \sum_{k = 1}^{n} g_k \nonumber \\ 
& \equiv &
\frac{p_{n + 1} -2}{n} \nonumber \\
& \sim & 
\ln (n)\,.    
\label{g-average}
\end{eqnarray}
Thus, we conclude that the average value of $\{g_n \}$ is about $\ln (n)$. 
A better average can be obtained by using better bounds of $p_k$.   
Since, $\left(\ln k + \ln \ln k \right) = \ln \left(k \ln k \right) < k$, 
we can also have an oversimplified statement that $\bar{g}_n < n$.  
  

\section{Statistical Approach to Andrica's Conjecture}
Let us consider the following sum using the Andrica's expression
\begin{eqnarray}
\sum_{k = 1}^{n} h_k 
& = & h_1 + h_2 + h_3 + \cdots + h_n \nonumber \\
& = & 
\left(\sqrt{p_2} - \sqrt{p_1} \right)
+ \left(\sqrt{p_3} - \sqrt{p_2} \right) 
+ \left(\sqrt{p_4} - \sqrt{p_3} \right)
+ \cdots \nonumber \\ 
& & \qquad \qquad \qquad 
+ \left(\sqrt{p_{n + 1}} - \sqrt{p_n} \right) \nonumber \\ 
& = & 
\sqrt{p_{n + 1}} - \sqrt{p_1}\,,  
\label{sum-h-exact}
\end{eqnarray}
where $\sqrt{p_1} = \sqrt{2}$.  We have used the telescoping property in the above summation.  
Next, we rewrite Eq.~(\ref{sum-h-exact}) as 
\begin{eqnarray}
\sum_{k = 1}^{n} h_k + \sqrt{2} = \sqrt{p_{n + 1}}\,.    
\label{sum-h-te-exact}
\end{eqnarray}
We divide the expression in Eq.~(\ref{sum-h-exact}) with $n$ and then rearrange to obtain 
the average value of the first $n$ gaps $\{ h \}$ denoted by $\bar{h}_n$  
\begin{eqnarray}
\bar{h}_n 
& \equiv & 
\frac{1}{n} \sum_{k = 1}^{n} h_k \nonumber \\ 
& \equiv & 
\frac{\sqrt{p_{n + 1}} - \sqrt{2}}{n}\,. 
\label{h-average-exact}
\end{eqnarray}
For the present, we shall use the oversimplification 
in Eq.~(\ref{upper-bound-of-pk-square}) to obtain 
\begin{eqnarray}
\bar{h}_n 
& \equiv &
\frac{\sqrt{p_{n + 1}} - \sqrt{2}}{n} \nonumber \\
& < & 
\frac{n + 1 - \sqrt{2}}{n} \nonumber \\ 
& = & 
1 - \frac{\sqrt{2} -1}{n} < 1\,.    
\label{h-average-simple}
\end{eqnarray}
Thus, we conclude that the average value of $\{h_n \}$ is less than unity.  
By construction $h_n$ is positive.  So, $\bar{h}_n \in (0\,, 1)$.  
This restriction on the average value of $\{h_n \}$ has implications.   
 
If any of the $h_k > 1$ then $\bar{h}_n \in (0\,, 1)$ ensures that 
there exists one or more complementary $h_{k^{\prime}} < 1$.  
For instance $h_k \in [1\,, 2)$ then $h_{k^{\prime}} < 1$ in order to sustain 
the $\bar{h}_n < 1$.  
We only know that the average value is less than unity.  Individual values may 
have lower or/and upper bounds!  
If $h_n$ happens to be bounded from below, then we shall require more than one 
$h_{k^{\prime}}$ to sustain the $\bar{h}_n \in (0\,, 1)$.  
If $h_k \in [2\,, 3)$, then we need at least two $h_{k^{\prime}} < 1$ to sustain  
the $\bar{h}_n < 1$.  Even if we include an $h_k \in [1\,, 2)$, we will still 
need at least two $h_{k^{\prime}} < 1$ in order to sustain the $\bar{h}_n \in (0\,, 1)$.  
Again, if the $h_n$ happens to be bounded from below, we shall require more than two 
$h_{k^{\prime}} < 1$ in order to sustain the average.  
Likewise, if $h_k \in [m\,, m + 1)$, then we require $m$ or more number of $h_{k^{\prime}} < 1$ 
due the restriction, $\bar{h}_n \in (0\,, 1)$.  
Thus, we can conclude that that half or more of the $\{h_n \}$ are less than one.  
This statistical inference is summarized in 
\begin{Theorem} 
{\bf (Andrica's Conjecture is satisfied by half or more prime pairs):} \label{theorem-andricas-conjecture-half}
Andrica's conjecture states that, for $p_n$ the $n$-th prime number, the inequality
$$
h_n \equiv \sqrt{p_{n + 1}} - \sqrt{p_n} < 1$$
holds for all $n$.  
Statistically, the average value of $\{ h_n \}$ is bounded as $\bar{h}_n \in (0\,, 1)$.  
Consequently, the inequality is satisfied by half or more of such prime pairs.    
\end{Theorem}

The statistical count can be improved to more than half, using standard statistical techniques. 


\section{Proof of Andrica's Conjecture}
The interrelated prime gap $g_n \equiv p_{n + 1} - p_n$ 
and the $h_n = \sqrt{p_{n + 1}} - \sqrt{p_n}$ 
can be analyzed using the average $\bar{h}_n$.  
We note the identity 
\begin{eqnarray}
\left(p_{n + 1} - p_n \right) 
= 
\left(\sqrt{p_{n + 1}} - \sqrt{p_n} \right) \left(\sqrt{p_{n + 1}} + \sqrt{p_n} \right)\,.  
\label{g-h-identity} 
\end{eqnarray}
Using the identity in Eq.~(\ref{g-h-identity}), we note the following 
\begin{eqnarray}
g_n = h_n \left(\sqrt{p_{n + 1}} + \sqrt{p_n} \right)\,.    
\label{g-h-n} 
\end{eqnarray}
Using the definition of averages and the relation in Eq.~(\ref{sum-h-te-exact}), 
we obtain the exact expressions 
\begin{eqnarray}
n \bar{h}_n 
& = & 
\sqrt{p_{n + 1}} - \sqrt{2} \nonumber \\ 
(n - 1) \bar{h}_{n - 1}  
& = & 
\sqrt{p_n} - \sqrt{2}\,,  
\label{h-bar-expressions} 
\end{eqnarray}
By construction in Eq.~(\ref{h-bar-expressions}), we have 
\begin{eqnarray}
n \bar{h}_n > (n - 1) \bar{h}_{n - 1}\,.  
\label{h-bar-inequality-basic} 
\end{eqnarray}
Then, we have 
\begin{eqnarray}
\bar{h}_n < \bar{h}_{n - 1}\,.  
\label{h-bar-inequality} 
\end{eqnarray}
It essentially means that $\bar{h}_n$ is a decreasing function of $n$.  
This can also be inferred from the fact that the average is bounded in the 
unit interval.  
Alternately, we can have an explicit function of $\bar{h}_n$ in terms of $n$.  
The required function is 
\begin{eqnarray}
\bar{h}_n 
& \equiv & 
\frac{1}{n} \sum_{k = 1}^{n} h_k \nonumber \\ 
& = &  
\frac{1}{n} \left\{\sqrt{p_{n + 1}} - \sqrt{2} \right\} \nonumber \\ 
& \sim & 
\frac{\sqrt{(n + 1) \ln (n + 1)} - \sqrt{2}}{n} \nonumber \\ 
& \sim & 
\frac{1}{\sqrt{\pi (n)}} < 1\,, 
\label{h-bar-structure}
\end{eqnarray}
where $\pi (n) = n/{\ln n}$ is the prime counting function.  
The average for large $n$ tends to zero. 
Again, we note that $\bar{h}_n$ is a decreasing function of $n$ with the 
restriction $\bar{h}_n \in (0\,, 1)$.  
From Eq.~(\ref{h-bar-inequality}) and Eq.~(\ref{h-bar-expressions}), 
we have 
\begin{eqnarray}
\bar{h}_n 
& < & 
\bar{h}_{n - 1} \nonumber \\ 
\frac{\sqrt{p_{n + 1}} - \sqrt{2}}{n} 
& < & 
\frac{\sqrt{p_n} - \sqrt{2}}{n - 1}\,.    
\label{hn-expressions-inequality} 
\end{eqnarray}
After some straightforward algebra on Eq.~(\ref{hn-expressions-inequality}), we have 
\begin{eqnarray} 
h_n 
& \equiv & 
\sqrt{p_{n + 1}} - \sqrt{p_n} \nonumber \\ 
& < & 
\frac{\sqrt{p_{n + 1}} - \sqrt{2}}{n} \nonumber \\ 
& = & 
\bar{h} \sim \frac{1}{\sqrt{\pi (n)}} < 1\,.    
\label{hn-less-than-one} 
\end{eqnarray}
Thus, we have the proof of the Andrica's conjecture.  
Moreover, we have been able to obtain a function for the Andrica's expression.  
This enabled us to note that $\lim h_n = 0$.  
Hence, the conjecture is proved completely.

\section{A Generalization of Andrica's Conjecture}
Let us consider the following generalization of Andrica's conjecture 
\begin{eqnarray}
h^{(x)}_n \equiv p^x_{n + 1} - p^x_n < 1\,.  
\label{andrica-generalized} 
\end{eqnarray}
We are considering a general power $x$ in place of the $1/2$.  
The range of $x$ shall also be determined during the proof.  
Then, our approach leads to the inequality 
\begin{eqnarray} 
h^{(x)}_n 
& \equiv 
& p^x_{n + 1} - p^x_n \nonumber \\ 
& < & 
\frac{p^x_{n + 1} - 2^x}{n} \nonumber \\ 
& \sim & 
\frac{p^x_{n + 1}}{n} < 1\,.    
\label{hn-less-than-one-andrica-generalized} 
\end{eqnarray}
This can be solved as 
\begin{eqnarray}
p^x_{n + 1} 
& < & 
n \nonumber \\ 
(n \ln n)^x 
& < & 
n \nonumber \\ 
\ln n
& < & 
n^{\frac{1}{x} - 1}\,. 
\label{inequality-for-andrica-generalized} 
\end{eqnarray}
The quantity $(\ln n)$ can always be suppressed by $n^b$ with $b > 0$ for $n \ge n_0$.  
In our case, $b = {\frac{1}{x} - 1}$.  
That is 
\begin{eqnarray}
\ln n 
& < & 
n^b\,, \quad \quad {\rm for ~} b > 0 \nonumber \\ 
n 
& < & 
e^{{n^b}}\,, \quad \quad {\rm for ~} n \ge n_0\,. 
\label{log-for-andrica-generalized}  
\end{eqnarray}
Consequently, $x \in (0\,, 1)$.  
For all $x$ outside this interval, the inequality 
in Eq.~(\ref{inequality-for-andrica-generalized}) fails as we are seeking solutions in $n$.   
The value of $n_0$ depends on the value by $b$. 
For $x \in (0\,, 1/2]$, $n_0 = 1$, since $b \ge 1$.  
In the remaining half of the interval $x \in (1/2\,, 1)$, $n_0 \ge 1$, since $b < 1$.
Hence, the generalized Andrica's conjecture is satisfied through the 
solvable inequality in~(\ref{inequality-for-andrica-generalized}).  
For $x = 1/2$, we have $\ln n < n$.  
This is the same oversimplification, which we have used in Eq.~(\ref{h-average-simple}).

\section{Concluding Remarks} 
We derived the average gap and showed that it is bounded in the interval $(0\,, 1)$. 
This strengthened the Andrica's conjecture and also affirmed that it is satisfied in 
half or more cases.  
By deriving an inequality for the averages, we are able to have a proof for the 
Andrica's conjecture.  
We also considered one of the generalization of Andrica's conjecture by relaxing 
the power in Andrica's expression from $1/2$ to a general $x \in (0\,, 1)$.  
This generalization was proved through an inequality, whose solution was also given in detail. 
In passing, we note that statistical techniques can be used to model the prime gaps.  
Of course, any such modelling has to be consistent with the established results.  

An immediate implication of Andrica's conjecture is on the 
prime gaps $g_n \equiv \left(p_{n + 1} - p_n \right)$.  
If we use the simple expression, the number $1$ for the Andrica's expression 
then, $g_n \sim 2 \sqrt{p_n} \sim 2 \sqrt{n \ln n}$ from Eq.~(\ref{g-h-n}).  
But now, we have a function for the Andrica's expression 
in Eq.~(\ref{hn-less-than-one}) leading to the bound  
\begin{eqnarray}
g_n 
& \equiv & 
p_{n + 1} - p_n \nonumber \\ 
& = &
h_n \left(\sqrt{p_{n + 1}} + \sqrt{p_n} \right) \nonumber \\
& < & 
2 h_n \sqrt{p_{n + 1}} \nonumber \\
& < & 
2 \left\{\frac{\sqrt{p_{n + 1}} - \sqrt{2}}{n} \right\} \sqrt{p_{n + 1}} \nonumber \\
& \sim & 
2 \frac{p_{n + 1}}{n} 
\sim 2 \frac{p_n}{n} \nonumber \\
& \sim & 2 \ln (n)\,.   
\label{g-bounds-andrica} 
\end{eqnarray}
This is twice the average value $\bar{g}_n \sim \ln (n)$ in Eq.~(\ref{g-average}). 
 
\bigskip
\bigskip

\noindent
{\bf Acknowledgements}: \\
My immediate interest in Andrica's conjecture was prompted by 
the masterpiece, {\em The M$\alpha$TH $\beta$OOK} by Clifford Alan Pickover 
(see pp 482-483 in \cite{MathBook}).  I recently received this book as a gift from my 
brother Prof. Farooq Ahmed Khan (University of West Georgia, Carrollton, GA, USA).  
I am grateful to Prof. Ramaswamy Jagannathan 
(Institute of Mathematical Sciences, Chennai, India),  
for numerous discussions in number theory.   



\end{document}